\documentclass{article}

\usepackage{mymacros}
\usepackage{hyperref, mathtools}
\usepackage{capt-of}
\usepackage{caption}
\usepackage{tikz-cd}
\begin{document}
\title{Veech Groups and Triangulations of Half-Dilation Pillowcases}
\author{Taro Shima}
\maketitle

\begin{abstract}
In this paper we consider the symmetries of triangulable half-dilation structures on the sphere with four singularities. We show that all such surfaces can be produced by a tetrahedral construction. Using this construction, we calculate each such surface's symmetry group in $\PSL(2, \R)$ called the Veech group.
\end{abstract}

\section{Introduction}

Similar to translation surfaces, there is a natural $\PSL(2, \R)$-action on the space of triangulable half-dilation structures. The Veech group of a triangulable half-dilation structure $P$ is $\PSL(P)$, the stabilizer of $P$ under the $\PSL(2, \R)$-action.\nextpar

Furthermore, $\PSL(2, \R)$ is the orientation-preserving isometry group of the hyperbolic plane $\mathbb{H}^2$. The Veech group $\PSL(P)$ can then be understood as a subgroup of orientation preserving isometries. In this paper, we prove the following result:

\begin{thm}
Given $P$, a triangulable half-dilation structure on the four-punctured sphere, its Veech Group $\PSL(P) < \PSL(2, \R)$ is such that $\mathbb{H}^2/\PSL(P)$ is a complete hyperbolic structure on the three-punctured sphere. Furthermore, every complete hyperbolic structure on the three-punctured sphere is represented.
\end{thm}

After introducing definitions, we proceed to show that every triangulable half-dilation structure on $S_{0,4}$ has a tetrahedral triangulation with the property that the sum of angles opposite every edge is less than $\pi$. We do this by using Delaunay triangulations. These triangulations can be represented as a triangle with gluings, and using this polygonal construction of the triangulable half-dilation structure, we calculate the Veech Group.

\section{Definitions}
Half-dilation surfaces are a certain generalization of half-translation surfaces, and can be constructed using half-dilations.
\begin{defn}[Half-Dilation]
A map $h:\C\to\C$ is called a half-dilation if it is of the form $z \mapsto az + b$ for $a\in\R_{\neq 0}$ and $b \in \C$.
\end{defn}

With half-translation surfaces, we obtain surfaces that contain cone points locally modeled by $k\pi$-cones. Half-dilation surfaces on the other hand have more information associated to them.

\begin{defn}[$k\pi$-Half-Dilation Cone]
A half-dilation $k\pi$-cone is defined by taking $k$ copies of the upper half-plane $\mathbb{U}_i$ and gluing $(-\infty, 0] \subset \mathbb{U}_n$ to $[0, \infty) \subset \mathbb{U}_{n+1 \mod k}$ via a half-dilation $z \mapsto a_iz$ where $a_i < 0$.
\end{defn}

We can now define a half-dilation structures -- Instead of requiring transition maps to be just translations or reflections we also allow for dilations. 

\begin{defn}[Half-Dilation Structure]
A half-dilation structure on a surface $X$ with a finite collection of singularities $\Sigma$ is an atlas on $X\wout\Sigma$ such that its transition maps are half-dilations. Furthermore, at each point $p \in \Sigma$, there is some $k$ and a homeomorphism of a neighborhood of $p$ to a neighborhood of the origin in a $k\pi$-half-dilation cone that is a half-dilation in local coordinates away from $p$.
\end{defn}

Though this allows for a much more diverse collection of surfaces, this generalization makes it more difficult to study. Some properties are preserved, for half-dilation surfaces are locally flat away from singularities, and these singularities respect the Gauss-Bonnet formula.

\begin{thm}[ {\cite[Theorem 1.1]{flat-surfaces}} ]
Let $\Sigma$ be the set of all cone singularities on a closed half-dilation surface $X$, with genus $g$. Then, for $p \in S$, let $k_p\pi$ be angle of the cone singularity. Then,
\[
2\pi(2 - 2g) = \sum_{p \in \sigma} 2\pi - k_p\pi
\]
\end{thm}

However, some properties aren't preserved. For example, translation surfaces have a natural metric associated to them. This construction unfortunately does not apply in the case of dilation surfaces. These surfaces have not yet been studied extensively, and understanding the Veech group will help us understand these surfaces better.\nextpar

\section{Triangulations and the Flipping Algorithm}

Triangulations are a very powerful tool to understand Veech group of half-dilation surfaces, as they will allow us to more concretely work with these surfaces.

\begin{defn}[Triangulation]
Given a surface $(X, h)$ with a half-dilation structure, a triangulation $\Delta_{(X,h)}$ of $(X,h)$ is a collection of finite collection of subsets $T_i \subset X$ such that:
\begin{align*}
1.& \textrm{the interiors of $T_i$ are pairwise disjoint and $X = \cup_i T_i$,}\\
2.& \textrm{the interior of each triangle admits a bijective local half-dilation to a triangle in $\R^2$,} \\
3.& \textrm{all vertices are points of singularity,}\\
4.& \textrm{the collection of points of singularity do not meet a triangle $T_i$ off its vertices.}
\end{align*}
\end{defn}

The existence of such a triangulation is proven in \cite{delaunay-partitions}. We will also be using a special triangulation called a Delaunay triangulation. It is special in that it has a desirable property regarding angles opposite edges in the triangulation.

\begin{defn}[Locally Delaunay]
An edge in a triangulation is incident to two triangles, and therefore, each edge has two angles it is opposite to it. An edge is called locally Delaunay if the sum of the angles opposite an edge is $\leq \pi$.
\end{defn}

\begin{defn}[Delaunay Triangulation]
A triangulation is Delaunay if all its edges are locally Delaunay. 
\end{defn}

\begin{thm}[{\cite[Theorem 10.8]{delaunay-partitions}}]
Every triangulable half-dilation surface structure on the 4-punctured sphere has a Delaunay triangulation.
\end{thm}
\begin{rmk}
The reader should be aware that Veech's Theorem 10.8 contains an error. The result he states does not hold in the generality claimed. But the theorem does hold in the context of triangulable half-dilation surfaces. See the discussion in \cite{dilation-veech}, in particular, appendix A.2.
\end{rmk}

There is also an algorithm to obtain a Delaunay triangulation from an arbitrary one through edge flips.
\begin{defn}
Two triangles in a triangulation that share an edge $e$ make a quadrilateral that has $e$ as a diagonal. Let the other diagonal be $e'$. Flipping the edge $e$ is defined by deleting the edge $e$ and drawing edge $e'$.
\end{defn}
The flipping algorithm is as follows:
\begin{align*}
1.&\ \textrm{Start with a triangulation $\Delta$.}\\
2.&\ \textrm{If all edges of the triangulation are Delaunay, stop.}\\
&\ \textrm{If not, there is an edge $e$ that is not locally Delaunay. Flip edge $e$.}\\
3.&\ \textrm{Go to Step 2.}
\end{align*}
Proving that this algorithm terminates in finite steps has been done for cone surfaces without dilational singularities in \cite{bobenko} and \cite{indermitte}, but does not generalize to half-dilation surfaces. We give a proof of finite termination that extends to similarity surfaces. In order to do this, we construct a functional on triangulations called the harmonic index.
\begin{defn}[Harmonic Index]
Let $P \subset \R^2$, be a polygon with sides $S_p = \set{s_1, \ldots, s_n}$ then the harmonic index $\operatorname{hrm}(P)$ is defined as,
\[
\operatorname{hrm}(P) = \sum_{s \in S_p} \frac{l(s)^2}{A(P)}
\]
where $l(s)$ is the length of side $s$ and $A(P)$ is the area of the polygon.
\end{defn}
The harmonic index of a triangulation $\Delta$ is defined by taking the sum of the harmonic indices of each triangle in the triangulation,
\[
\operatorname{hrm}(\Delta) = \sum_{T_i \in \Delta} \operatorname{hrm}(T_i)
\]
There are two important properties of the harmonic index: 
\begin{prop}[\cite{musin}]\hypertarget{prop3.8}
Given a non-Delaunay triangulation $\Delta$, and $\Delta'$, a triangulation resulting from flipping an edge of $\Delta$ that is not locally Delaunay, 
\[
\operatorname{hrm}(\Delta') \leq \operatorname{hrm}(\Delta)
\]
\end{prop}
We extend the notion of the harmonic index to the collection of triangulations of similarity surfaces. \nextpar
In {\cite[Definition 6.11]{flat-surfaces}}, Veech defines $\cal{S}_1(g,n)$, the collection of all similarity structures $M_{g,n}$ up to isotopy with non-negative cone angles. By {\cite[Theorem 1.16]{flat-surfaces}}, , it is a complex analytic manifold. Veech describes local homeomorphisms from $\mathbb{H}^n$ to $\cal{S}_1(g,n)$. We follow Veech's construction in Section 5. In order to describe this homeomorphism, we first fix a combinatorial triangulation $t$ with $n$ triangles. The collection of triangulations $\Omega(t)$ is an open set of $\mathbb{H}^n$ where each $(z_1, \ldots, z_n) \in \Omega(t)$ is associated to the collection of $n$ triangles connecting $-1, 1, z_i  \in \mathbb{H}$ as shown in Figure 1, with the appropriate gluings given by the combinatorial triangulation.

\begin{center}
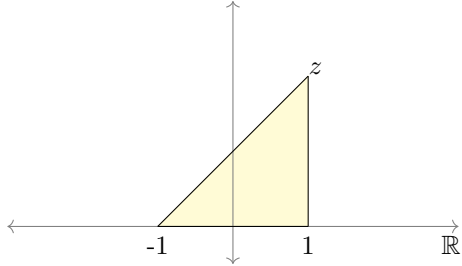

\begin{tikzpicture}
\filldraw[yellow!20!white] (-1,0) -- (1,2) -- (1,0) -- cycle;
\draw[black!50!white, <->] (-3, 0) -- (3,0);
\draw[black!50!white, <->] (0,3) -- (0,-0.5);

\draw[-] (-1,0) -- (1,2);
\draw[-] (1,0) -- (1,2);
\draw[-] (-1,0) -- (1,0);

\draw (2.9, -0.25) node {$\R$};
\draw (1, -0.25) node {1};
\draw (-1, -0.25) node {-1};
\draw (1.1, 2.1) node {$z$};
\end{tikzpicture}
\captionof{figure}{Triangle associated to $z \in \mathbb{H}$}
\end{center}
Veech defines a map $\sigma_t : \Omega(t) \to \cal{S}_1(g,n)$ , which by Theorem 1.16, is a homeomorphism onto its image.

Thus, fixing a combinatorial triangulation, we can define a continuous map on subsets of $\cal{S}_1(g,n)$,
\[
\operatorname{hrm}\circ\ \sigma_t^{-1}: \sigma_t(\Omega(t)) \to \R.
\]
We will use this functional in order to show that the flipping algorithm terminates in finitely many steps.
\begin{thm}
The flipping algorithm terminates in finitely many steps.
\end{thm}
\begin{proof}
We proceed with a proof by contradiction. Assume not.\nextpar
Then, for a triangulation $\Delta_0$ of a triangulable similarity structure on a surface $X$, there is some infinite sequence of flips of non-Delaunay edges $\set{f_1, f_2, \ldots}$ that result in triangulations $\set{f_1(\Delta_0), f_2f_1(\Delta_0), \ldots}$. Let $f_n\cdots f_1 = g_n$ and $g_n(\Delta_0) = \Delta_{n+1}$. There are only finitely many triangulations up to combinatorial type, as there a finite number of vertices so we can pass to a subsequence of triangulations $\set{\Delta'_i}$ so that each triangulation in $\set{\Delta_i'}$ are of the same combinatorial type. Let this combinatorial triangulation of $X$ be $t$. \nextpar
Now we show that $\set{\Delta_i'}$ has an accumulation point in $\Omega(t)$. Consider the set $A = \set{\mathbb{H} \in \Omega(t) \st \operatorname{hrm}(\Delta) \leq \operatorname{hrm}(\Delta_0)}$. By \hyperlink{prop3.8}{Proposition 3.8}, the sequence $\set{\Delta_i'} \subset A$. Furthermore, $A \subset \set{z \in \mathbb{H} \st \operatorname{hrm}(z) \leq \operatorname{hrm}(\Delta_0)}^n \subset \mathbb{H}^n$, since the harmonic index of each of the triangles cannot exceed the harmonic index of the triangulation. Furthermore, $\set{z \in \Omega(t) \st \operatorname{hrm}(z) \leq \operatorname{hrm}(\Delta_0)}$ is a compact set, and since the finite product of compact sets is compact, we conclude that $\set{\Delta_i'}$ must be a sequence contained in a compact set. We can then pass to a subsequence $\set{\Delta_{i_k}'}$ that converges to some $\Delta_T \in \Omega(t)$. \nextpar
Using the homemomorphism $\sigma_t$, we obtain a converging sequence $\set{\sigma_t(\Delta_{i_k}')} \subset \cal{S}_1(g,n)$. Since each $\Delta_{i_k}'$ triangulates the same similarity structure, $\sigma_t(\Delta_{i_k}')$ and $\sigma_t(\Delta_{i_j}')$ are equivalent up to an element of the mapping class group.
\nextpar
However, this contradicts {\cite[Theorem 7.22]{flat-surfaces}}, as the mapping class group acts discretely on $\mathcal{S}_1(g,n)$. We conclude that the algorithm terminates in a finite number of steps.
\end{proof}

\section{Triangulable Half-Dilation Pillowcases}

Flat pillowcases are the one of the first examples of half-translation surfaces that one might think of. It's polygonal gluing looks like so:

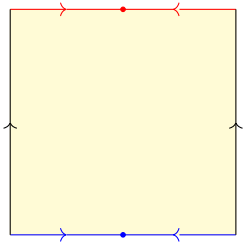
\begin{center}
\begin{tikzpicture}
\filldraw[yellow!20!white] (-1.5, -1.5) rectangle (1.5,1.5);
\draw[red, ->] (-1.5, 1.5) -- (-0.75,1.5);
\draw[red, -] (-0.75, 1.5) -- (0,1.5);

\draw[red, -<] (0, 1.5) -- (0.75,1.5);
\draw[red, -] (0.75, 1.5) -- (1.5,1.5);

\draw[blue, ->] (-1.5, -1.5) -- (-0.75,-1.5);
\draw[blue, -] (-0.75, -1.5) -- (0,-1.5);

\draw[blue, -<] (0, -1.5) -- (0.75,-1.5);
\draw[blue, -] (0.75, -1.5) -- (1.5,-1.5);

\draw[->] (-1.5, -1.5) -- (-1.5, 0);
\draw[-] ( -1.5, 0) -- (-1.5,1.5);

\draw[->] (1.5, -1.5) -- (1.5, 0);
\draw[-] (1.5, 0) -- (1.5,1.5);

\filldraw[blue] (0, -1.5) circle (0.03);

\filldraw[red] (0, 1.5) circle (0.03);

\end{tikzpicture}
\end{center}

We generalize this into a half-dilation surface by considering a surface with 4 half-dilation singularities. We proceed to understand what a triangulation of a triangulable half-dilation pillowcase may look like.

\begin{prop}
A triangulation of a triangulable half-dilation pillowcase cannot contain an edge that connects a vertex $p$ to itself.
\end{prop}
\begin{proof}
Without loss of generality, after cutting along the loop and flattening it on to $\R^2$, we have the following:
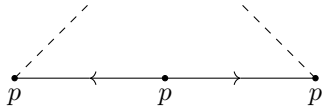
\begin{center}
\begin{tikzpicture}
\draw[->] (0,0) -- (1,0);
\draw (1,0) -- (2,0);
\draw[->] (0,0) -- (-1,0);
\draw (-1,0) -- (-2,0);

\filldraw (0,0) circle (1pt);
\filldraw (-2,0) circle (1pt);
\filldraw (2,0) circle (1pt);

\draw (0, -0.25) node {$p$};
\draw (2, -0.25) node {$p$};
\draw (-2, -0.25) node {$p$};

\draw[dashed] (2,0) -- (1, 1);
\draw[dashed] (-2, 0) -- (-1, 1);
\end{tikzpicture}
\end{center}

However, because the interior of each subset must be isometric to a triangle, regardless of which other vertices $p$ is connected to, the angle around $p$ will be more than $\pi$. This is not allowed in our triangulation, as all singularities are $\pi$-singularities.
\end{proof}

\begin{prop}
A triangulable half-dilation pillowcase admits at most two types of minimal triangulations.
\end{prop}
\begin{proof}
By definition, there are 4 vertices in the triangulation. Furthermore, there must be 4 triangles in the triangulation since very triangle contributes $\pi$ angle and the sum of the angle around each vertex must be $4\pi$. \nextpar

Let ${v_1, v_2, v_3, v_4}$ be the set of vertices. If we let $f(v_i)$ be the number of faces around the vertex $v_i$,
\[
\sum_{i = 1}^4 f(v_i) = 12,
\]
since each triangular face is triple counted. Therefore, we can focus on partitions of 12 into 4 integers,
\[
\set{f(v_1), f(v_2), f(v_3), f(v_4)}.
\] 
However, we note that none of the integers can be 1, since it would imply that there is a triangle in the triangulation that has an angle of $\pi$. We have the additional constraint that a partition cannot have an integer greater than 4, for if a vertex had more than 4 faces around it, it would imply that an edge connects a vertex to itself.
So, we reduce to the following partitions:
\[
\set{2, 2, 4, 4},\quad \set{2, 3, 3, 4},\quad \set{3, 3, 3, 3}.
\]

We now separate into case work by considering the number of edges between two vertices. If two vertices are connected by $n\geq 3$ edges, topologically, the edges split the sphere into $n$ discs that have 2 points of singularity on the boundary. Regardless of where the other singularities lie, there will be one disc that has no singularities in its interior. This results in a triangulation with a 2-gon, which is a contradiction. \nextpar

Now, if two vertices are connected by 2 edges, the edges split the sphere into 2 discs that have 2 points of singularity on the boundary. In order to not create a 2-gon, there must be exactly 1 singularity on the interior of each disc. Since an edge cannot connect a vertex to itself, the singularity on the interior must be connected by an edge to each of the singularities on the boundary. This results in a triangulation with the partition $\set{2, 2, 4, 4}$. Such a triangulation must look like so:

\begin{center}
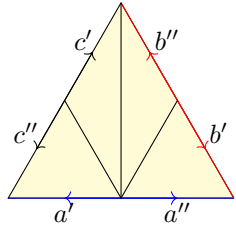

\begin{tikzpicture} [scale = 0.75]

\filldraw[yellow!20!white] (({2},0) -- (0,{2*sqrt(3)}) -- ({-2},0) -- ({2},0) -- cycle;
\draw[-] ({2},0) -- (0,{2*sqrt(3)}) -- ({-2},0) -- ({2},0);

\draw[blue] ({2},0) -- (-2,0);
\draw[blue, ->] (0,0) -- (1,0);
\draw[blue, ->] (0,0) -- (-1,0);

\draw[red, ->] ({1},{sqrt(3)}) -- (1.5, {0.5*sqrt(3)});
\draw[red, ->] ({1},{sqrt(3)}) -- (0.5, {1.5*sqrt(3)});
\draw[red] ({2},0) -- (0,{2*sqrt(3)});

\draw[->] ({-1},{sqrt(3)}) -- (-1.5, {0.5*sqrt(3)});
\draw[->] ({-1},{sqrt(3)}) -- (-0.5, {1.5*sqrt(3)});

\draw (1, -0.25) node {$a''$};
\draw (-1, -0.25) node {$a'$};

\draw (2/3+0.15, {(4/3)*sqrt(3) + 0.5}) node {$b''$};
\draw (4/3 + 0.4, {(2/3)*sqrt(3)}) node {$b'$};

\draw (-2/3, {(4/3)*sqrt(3) + 0.5}) node {$c'$};
\draw (-4/3 - 0.35, {(2/3)*sqrt(3)}) node {$c''$};

\draw[-] (-1,{sqrt(3)}) -- (0,0) -- (1,{sqrt(3)});
\draw[-] (0,0) -- (0,{2*sqrt(3)});

\end{tikzpicture}
\captionof{figure}{\textrm{$\set{2,3,3,4}$-Triangulation of a Half-Dilation Structure on $S_{0,4}$}}
\end{center}
If there are no two vertices that have more than 1 edge between them, the only triangulation that satisfies this condition is the tetrahedral triangulation with the partition $\set{3, 3, 3, 3}$, like so:

\begin{center}
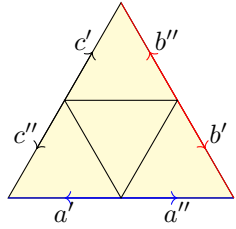

\begin{tikzpicture}[scale = 0.75]

\filldraw[yellow!20!white] (({2},0) -- (0,{2*sqrt(3)}) -- ({-2},0) -- ({2},0) -- cycle;
\draw[-] ({2},0) -- (0,{2*sqrt(3)}) -- ({-2},0) -- ({2},0);

\draw[blue] ({2},0) -- (-2,0);
\draw[blue, ->] (0,0) -- (1,0);
\draw[blue, ->] (0,0) -- (-1,0);

\draw[red, ->] ({1},{sqrt(3)}) -- (1.5, {0.5*sqrt(3)});
\draw[red, ->] ({1},{sqrt(3)}) -- (0.5, {1.5*sqrt(3)});
\draw[red] ({2},0) -- (0,{2*sqrt(3)});

\draw[->] ({-1},{sqrt(3)}) -- (-1.5, {0.5*sqrt(3)});
\draw[->] ({-1},{sqrt(3)}) -- (-0.5, {1.5*sqrt(3)});

\draw (1, -0.25) node {$a''$};
\draw (-1, -0.25) node {$a'$};

\draw (2/3+0.15, {(4/3)*sqrt(3) + 0.5}) node {$b''$};
\draw (4/3 + 0.4, {(2/3)*sqrt(3)}) node {$b'$};

\draw (-2/3, {(4/3)*sqrt(3) + 0.5}) node {$c'$};
\draw (-4/3 - 0.35, {(2/3)*sqrt(3)}) node {$c''$};

\draw[-] ({1},{sqrt(3)}) -- (-1,{sqrt(3)}) -- (0,0) -- (1,{sqrt(3)});

\end{tikzpicture}
\captionof{figure}{\textrm{Tetrahedral Triangulation of a Half-Dilation Structure on $S_{0,4}$}}
\label{tet_trig}
\end{center}
Any other triangulation must have the property that no two vertices have more than 1 edge between them, and also that no two vertices have exactly 1 edge between them. Such a triangulation cannot exist, so we conclude that a $\set{2,3,3,4}$-triangulation cannot exist.
\end{proof}
Although this classifies triangulations, these are not necessarily Delaunay triangulations. We introduce a proposition that determines when a Delaunay triangulation is unique.

\begin{prop}
The following Delaunay triangulation is unique when \[
\alpha + \alpha' < \pi,\quad\beta + \beta' < \pi,\quad \gamma + \gamma' < \pi
\]
where the angles are depicted as in Figure 4
\begin{center}
\begin{tikzpicture}[scale=0.95]

\filldraw[yellow!20!white] (({2},0) -- (0,{2*sqrt(3)}) -- ({-2},0) -- ({2},0) -- cycle;
\draw[-] ({2},0) -- (0,{2*sqrt(3)}) -- ({-2},0) -- ({2},0);

\draw[blue] ({2},0) -- (-2,0);
\draw[blue, ->] (0,0) -- (1,0);
\draw[blue, ->] (0,0) -- (-1,0);

\draw[red, ->] ({1},{sqrt(3)}) -- (1.5, {0.5*sqrt(3)});
\draw[red, ->] ({1},{sqrt(3)}) -- (0.5, {1.5*sqrt(3)});
\draw[red] ({2},0) -- (0,{2*sqrt(3)});

\draw[->] ({-1},{sqrt(3)}) -- (-1.5, {0.5*sqrt(3)});
\draw[->] ({-1},{sqrt(3)}) -- (-0.5, {1.5*sqrt(3)});

\draw (1, -0.25) node {$a''$};
\draw (-1, -0.25) node {$a'$};

\draw (2/3+0.15, {(4/3)*sqrt(3) + 0.5}) node {$b''$};
\draw (4/3 + 0.4, {(2/3)*sqrt(3)}) node {$b'$};

\draw (-2/3, {(4/3)*sqrt(3) + 0.5}) node {$c'$};
\draw (-4/3 - 0.35, {(2/3)*sqrt(3)}) node {$c''$};

\draw[-] ({1},{sqrt(3)}) -- (-1,{sqrt(3)}) -- (0,0) -- (1,{sqrt(3)});

\draw (-1.5,0) arc (0:60:0.5);
\draw (-1.35, 0.35) node {$\alpha$};
\draw (1.5, 0) arc (180:120:0.5);
\draw (1.35, 0.35) node {$\beta$};

\draw (-0.25, {1.75*sqrt(3)}) arc (240:300:0.5);
\draw (0, {1.75*sqrt(3)-0.35}) node {$\gamma$};

\draw (({-1+0.35},{sqrt(3)}) arc (360:300:0.35);
\draw (({1-0.35},{sqrt(3)}) arc (180:240:0.35);
\draw (({0.35*0.5},{0.5*0.35*sqrt(3)}) arc (60:120:0.35);

\draw (0, 0.6) node {$\gamma'$};
\draw (0 + 0.45 , {0.75 + 0.45*sqrt(3) - 0.05}) node {$\alpha'$};
\draw ({-0.45}, {0.75 + 0.45*sqrt(3) - 0.05}) node {$\beta'$};
\end{tikzpicture}
\captionof{figure}{}
\end{center}

for $a = \tfrac{a''}{a'}, b = \tfrac{b''}{b'}, c = \tfrac{c''}{c'} \in \R_{>0}$ and $\alpha, \beta \in \R_{>0}$ such that $0
<  \alpha + \beta < \pi$.

\end{prop}
\begin{proof}

First, we show that a $\set{2, 2, 4, 4}$-triangulation cannot be unique. Consider the following:
\begin{center}
\begin{tikzpicture}[scale=0.75]

\filldraw[yellow!20!white] (({2},0) -- (0,{2*sqrt(3)}) -- ({-2},0) -- ({2},0) -- cycle;
\draw[-] ({2},0) -- (0,{2*sqrt(3)}) -- ({-2},0) -- ({2},0);

\draw[blue] ({2},0) -- (-2,0);
\draw[blue, ->] (0,0) -- (1,0);
\draw[blue, ->] (0,0) -- (-1,0);

\draw[red, ->] ({1},{sqrt(3)}) -- (1.5, {0.5*sqrt(3)});
\draw[red, ->] ({1},{sqrt(3)}) -- (0.5, {1.5*sqrt(3)});
\draw[red] ({2},0) -- (0,{2*sqrt(3)});

\draw[->] ({-1},{sqrt(3)}) -- (-1.5, {0.5*sqrt(3)});
\draw[->] ({-1},{sqrt(3)}) -- (-0.5, {1.5*sqrt(3)});

\draw (1, -0.25) node {$a'$};
\draw (-1, -0.25) node {$a$};

\draw (2/3+0.15, {(4/3)*sqrt(3) + 0.5}) node {$b'$};
\draw (4/3 + 0.4, {(2/3)*sqrt(3)}) node {$b$};

\draw (-2/3, {(4/3)*sqrt(3) + 0.5}) node {$c$};
\draw (-4/3 - 0.35, {(2/3)*sqrt(3)}) node {$c'$};

\draw[-] (-1,{sqrt(3)}) -- (0,0) -- (1,{sqrt(3)});
\draw[-] (0,0) -- (0,{2*sqrt(3)});

\draw (1.25, {0.75*sqrt(3)}) arc (-60:-120:0.5);
\draw (-1.25, {0.75*sqrt(3)}) arc (-120:-60:0.5);

\draw ({1},{sqrt(3)-0.75}) node {$\phi$};
\draw ({-1},{sqrt(3)-0.75}) node {$\theta$};

\end{tikzpicture}
\captionof{figure}{}
\end{center}

A triangulation is Delaunay if the sum of angles opposite edges are all less than or equal to $\pi$. Therefore, a Delaunay triangulation is unique if the sum of angles opposite edges are all less than $\pi$. This cannot be the case in the $\set{2, 2, 4, 4}$-triangulation, since $\phi + \theta < \pi$ implies the sum of the supplemental angles are greater than $\pi$, and vice versa. Therefore, $\phi + \theta = \pi$. This implies that the triangulation can be flipped -- making the triangulation non-unique.

Now consider a tetrahedral triangulation. Instead of checking the sum of the angles opposite all 6 edges, we can check the sum of the angles opposite three edges, as seen in the proposition. This follows from Euclidean geometry. In particular, $\alpha + \alpha'$ is equal to the sum of the angle opposite the edges $b'$ and $b''$ since, they lie supplemental to the double ringed angles as in Figure 5.
\begin{center}
\begin{tikzpicture}[scale=0.75]

\filldraw[yellow!20!white] (({2},0) -- (0,{2*sqrt(3)}) -- ({-2},0) -- ({2},0) -- cycle;
\draw[-] ({2},0) -- (0,{2*sqrt(3)}) -- ({-2},0) -- ({2},0);

\draw[blue] ({2},0) -- (-2,0);
\draw[blue, ->] (0,0) -- (1,0);
\draw[blue, ->] (0,0) -- (-1,0);

\draw[red, ->] ({1},{sqrt(3)}) -- (1.5, {0.5*sqrt(3)});
\draw[red, ->] ({1},{sqrt(3)}) -- (0.5, {1.5*sqrt(3)});
\draw[red] ({2},0) -- (0,{2*sqrt(3)});

\draw[->] ({-1},{sqrt(3)}) -- (-1.5, {0.5*sqrt(3)});
\draw[->] ({-1},{sqrt(3)}) -- (-0.5, {1.5*sqrt(3)});

\draw (1, -0.25) node {$a''$};
\draw (-1, -0.25) node {$a'$};

\draw (2/3+0.15, {(4/3)*sqrt(3) + 0.5}) node {$b''$};
\draw (4/3 + 0.4, {(2/3)*sqrt(3)}) node {$b'$};

\draw (-2/3, {(4/3)*sqrt(3) + 0.5}) node {$c'$};
\draw (-4/3 - 0.35, {(2/3)*sqrt(3)}) node {$c''$};

\draw[-] ({1},{sqrt(3)}) -- (-1,{sqrt(3)}) -- (0,0) -- (1,{sqrt(3)});

\draw (-1.5,0) arc (0:60:0.5);
\draw (-1.35, 0.35) node {$\alpha$};

\draw (({1-0.35},{sqrt(3)}) arc (180:240:0.35);

\draw (({-1+0.35},{sqrt(3)}) arc (360:240:0.35);
\draw (({-1+0.3},{sqrt(3)}) arc (360:240:0.3);

\draw (({0.35*0.5},{0.5*0.35*sqrt(3)}) arc (60:180:0.35);
\draw (({0.3*0.5},{0.5*0.3*sqrt(3)}) arc (60:180:0.3);

\draw (0 + 0.45 , {0.75 + 0.45*sqrt(3) - 0.05}) node {$\alpha'$};

\end{tikzpicture}
\captionof{figure}{}
\end{center}

This argument can be repeated for $\alpha + \alpha'$ and $\beta + \beta'$, so we conclude that it suffices to check,
\[
\alpha + \alpha' < \pi, \quad \beta + \beta' < \pi,\quad \gamma + \gamma' < \pi
\]
to determine that a Delaunay triangulation is unique.

\end{proof}
We can further analyze when triangulations are not unique. Note that of the 3 inequalities:
\[
\alpha + \alpha' < \pi,\quad \beta + \beta' < \pi,\quad \gamma + \gamma' < \pi,
\]
no two can equality hold. This is because, if $\alpha + \alpha' = \pi$ and $\beta + \beta' = \pi$, we would end up with a quadrilateral with more than $2\pi$ angle, as shown in Figure 7.
\begin{center}
\begin{tikzpicture}[scale=0.75]

\filldraw[yellow!20!white] (({2},0) -- (0,{2*sqrt(3)}) -- ({-2},0) -- ({2},0) -- cycle;
\filldraw[yellow!60!white] (({2},0) -- ({1},{sqrt(3)}) -- ({-1},{sqrt(3)}) -- ({-2},0) -- ({2},0) -- cycle;
\draw[-] ({2},0) -- (0,{2*sqrt(3)}) -- ({-2},0) -- ({2},0);

\draw[blue] ({2},0) -- (-2,0);
\draw[blue, ->] (0,0) -- (1,0);
\draw[blue, ->] (0,0) -- (-1,0);

\draw[red, ->] ({1},{sqrt(3)}) -- (1.5, {0.5*sqrt(3)});
\draw[red, ->] ({1},{sqrt(3)}) -- (0.5, {1.5*sqrt(3)});
\draw[red] ({2},0) -- (0,{2*sqrt(3)});

\draw[->] ({-1},{sqrt(3)}) -- (-1.5, {0.5*sqrt(3)});
\draw[->] ({-1},{sqrt(3)}) -- (-0.5, {1.5*sqrt(3)});

\draw[-] ({1},{sqrt(3)}) -- (-1,{sqrt(3)}) -- (0,0) -- (1,{sqrt(3)});

\draw (-1.5,0) arc (0:60:0.5);
\draw (-1.35, 0.35) node {$\alpha$};
\draw (1.5, 0) arc (180:120:0.5);
\draw (1.35, 0.35) node {$\beta$};

\draw (({-1+0.35},{sqrt(3)}) arc (360:300:0.35);
\draw (({1-0.35},{sqrt(3)}) arc (180:240:0.35);

\draw (0 + 0.45 , {0.75 + 0.45*sqrt(3) - 0.05}) node {$\alpha'$};
\draw ({-0.45}, {0.75 + 0.45*sqrt(3) - 0.05}) node {$\beta'$};
\end{tikzpicture}
\captionof{figure}{}
\end{center}
So, we conclude the following.

\begin{prop}
If a triangulable half-dilation pillowcase has a non-unique Delaunay triangulation, exactly one of the following conditions hold:
\begin{align*}
&\alpha + \alpha' = \pi & &&&\alpha + \alpha' < \pi && &&\alpha + \alpha' < \pi \\
&\beta + \beta' < \pi && or &&\beta + \beta' = \pi & & or && \beta + \beta' < \pi \\
&\gamma + \gamma' < \pi &&&& \gamma + \gamma' < \pi &&&& \gamma + \gamma' = \pi
\end{align*}
Furthermore, such a triangulable half-dilation surface has 4 possible Delaunay triangulations.
\end{prop}
\begin{proof}
As per the discussion above, if any two pairs of angles sum to exactly $\pi$ angle, we would end up with a quadrilateral with more than $2\pi$ angle, which is a contradiction.\nextpar

Now assume that one pair of angles in the tetrahedral triangulation sums to $\pi$ degrees. Then, there are two edges for which the angles opposite the edge sum to $\pi$ angle. Therefore, there are 2 places that can be flipped to get another Delaunay triangulation. Thus, we end up with 4 possible triangulations for such a half-dilation pillowcase.\nextpar

In particular, there are 2 tetrahedral triangulations, where 0 edges are flipped or both edges are flipped, and 2 $\set{2, 2, 4, 4}$-triangulations, where only 1 of the edges are flipped.
\end{proof}

\section{Veech Groups of Half-Dilation Pillowcases}

In the study of translation surfaces understanding the stabilizer of $(X, \omega)$ under the $\SL(2, \R)$-action, also known as the Veech Group, offers a lot information about the translation surfaces in question. For example, Veech proved the following:
\begin{thm}[Veech Dichotomy]
Let $(X, \omega)$ be a translation surface. Suppose $\SL(X, \omega)$ is a lattice in $\SL(2, \R)$. Then for each direction $\theta$, the flow $F_\theta$ is either periodic or uniquely ergodic.
\end{thm}
We will calculate an analog for the Veech Group for triangulable half-dilation pillowcases. In particular, instead of finding the derivative of affine diffeomorphims, we will find the derivatives of $\SL(2,\R)$-diffeomorphisms.
\begin{defn}[Affine Automorphisms]
Let $q:\GL(2,\R) \to \PSL(2,\R)$ be the quotient map and $X$ a half-dilation surface. A diffeomorphism $f:X \to X$ is called a $\SL(2,\R)$ diffeomorphism if $q(D_xf)$ is constant throughout $X$ off of singularities
\end{defn}

\subsection{Marked Triangles}

We briefly go over marked triangles, which will be a way to put coordinates on the space of half-dilation pillowcases.

\begin{defn}[Marked Triangle]
Let $\Delta \subset \R^2$ be a Euclidean triangle whose boundary is a union of three line segments $l_i$ where we normalize $\abs{l_1} = 1$. A marked triangle is $(\Delta, p_1, p_2, p_3)$ where $p_i \in l_i\wout\partial l_i$
\end{defn}
Let $\Delta^*$ be the collection of all marked triangles modulo a translation and let $\scr{T} \subset (\R^2 \times \R)^3$ be the collection of triples of the form $\big((v_1, \lambda_1), (v_2, \lambda_2), (v_3, \lambda_3)\big)$ where $v_i = (x_1, y_i) \in \R^2,$ the parameters $ \lambda_i \in \R_{>0}$, the vectors $v_1, v_2$ are linearly independent, and $v_1 + v_2 + v_3 = 0$. We can specify a coordinate system on marked triangles by defining a map $\cal{T}:\scr{T}\to \Delta^*$. \nextpar

The map $\mathcal{T}$ sends a triple $(v_1, \lambda_1), (v_2, \lambda_2), (v_3, \lambda_3)$ to the triangle defined by connecting the origin, and the points with coordinates $v_1$ and $v_1 + v_2$. Using the real numbers $\lambda_1, \lambda_2, \lambda_3$ we can take the points of singularities $p_1, p_2, p_3$ to be $\tfrac{1}{\lambda_1 + 1}v_1, \tfrac{1}{\lambda_2 + 1}v_2 + v_1$ and $\tfrac{1}{\lambda_3 + 1}v_3 + v_2 + v_1$ for $\lambda_1, \lambda_2, \lambda_3$ respectively. Therefore the triple $((v_1, \lambda_1),(v_2, \lambda_2),(v_3, \lambda_3))$ is sent to the marked triangle as depicted in Figure 8.

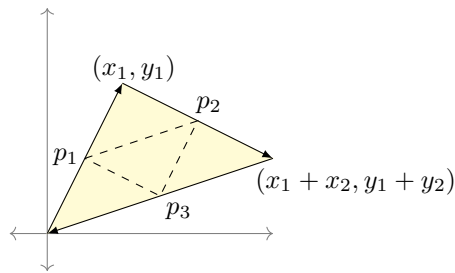
\begin{center}
\begin{tikzpicture}
\draw[black!50!white, <->] (-0.5, 0) -- (3,0);
\draw[black!50!white, <->] (0,3) -- (0,-0.5);

\filldraw[yellow!20!white] (0,0) -- (1,2) -- (3,1) -- cycle;

\draw[-latex] (0,0) -- (1,2);
\draw[-latex] (1,2) -- (3,1);
\draw[-latex] (3,1) -- (0,0);

\draw (1.15, 2.2) node {$(x_1, y_1)$};
\draw (4.10, 0.70) node {$(x_1 + x_2, y_1 + y_2)$};

\draw[dashed] (0.5,1) -- (2, 1.5) -- (1.5, 0.5) -- cycle;
\draw (0.5 - 0.25, {0.5*sqrt(3) + 0.15}) node {$p_1$};
\draw (1.5 + 0.65, {0.5*sqrt(3) + 0.85}) node {$p_2$};
\draw (1.75, 0.3) node {$p_3$};
\end{tikzpicture}
\captionof{figure}{Marked triangle associated to $\set{(v_1, \lambda_1), (v_2, \lambda_2), (v_3, \lambda_3)}$ in $\R^2$}
\end{center}

We consider a half-dilation structure on the pillowcase, with the singularities labeled by $\set{p_1,p_2,p_3,p_4}$. Suppose this surface has a triangulation homeomorphic to the tetrahedron. Then we can cut along the edges connecting to the $p_4$ vertex and flatten the result onto $\R^2$.  As the singularities have angle $\pi$, we obtain a triangle in $\R^2$. Without loss of generality, we can take a vertex of the edge containing the singularity $p_1$ to be the origin. We label the singularities $\set{p_1, p_2, p_3}$ on the edges of this triangle to obtain a marked triangle, as depicted in Figure 8. \nextpar

On the other hand by gluing this tetrahedron using the gluing scheme as shown in Figure 1, we obtain a $S_{0,4}$ with a half-dilation structure. Let $\Delta^*(S_{0,4})$ be the collection of all tetrahedrally triangulated half-dilation pillowcases. We endow $\Delta^*(S_{0,4})$ with the quotient topology, since $\cal{G}$ is surjective.
\begin{center}
\begin{tikzcd}
\scr{T} \ar[r,"\cal{T}"] & \Delta^* \ar[r, "\cal{G}"] & \Delta^*(S_{0,4})
\end{tikzcd}
\end{center}

\subsection{Half-Dilation Pillowcases}

We can take an equivalence relation on the collection of all tetrahedral triangulations of half-dilation pillowcases in order to obtain the collection of all half-dilation pillowcases.
\begin{prop}
Two tetrahedral triangulations $\Delta_1$ and $\Delta_2$ represent the same half-dilation surface if there is a sequence of edge flips and cut and paste operations to get from one tetrahedrally triangulated surface to another.
\end{prop}
\begin{proof}
By theorems 2.8 and 2.11, there is a sequence of finite edge flips to get from any tetrahedral triangulation to a Delaunay triangulation $\Delta_D$. By taking the sequence of flips from $\Delta_1$ to $\Delta_D$ and reversing the flips to get to $\Delta_D$ to $\Delta_2$, we obtain our desired sequence of flips.\nextpar

After taking edge flips on $\Delta_1$ to get to $\Delta_2$, to get a triangular representation of $\Delta_2$ in $\R^2$, we can cut along the edges resulting from the edge flips, and paste after scaling, rotating by $\pi$ to get to the triangular representation of $\Delta_2$.
\end{proof}

This proposition implies that in order to calculate the stabilizer of a half-dilation pillowcase under the $\PSL(2, \R)$-action, it suffices to understand when $X$ and $A\cdot X$ differ by a cut-scale-paste operation. We will answer this by determining affine maps associated to each cut-scale-paste operation.

\nextpar
We define a special map that takes edge flips on two disjoint pairs of triangles, and then cutting, scaling and pasting so that a tetrahedral triangulation is represented in $\R^2$. There are three such maps, $\Phi_1, \Phi_2, \Phi_3: \Delta^* \to \Delta^*$ as shown in Figure 9. We give a more formal definition of the maps $\Phi_i:\Delta^*\to \Delta^*$. Given a marked triangle associated to a tetrahedral triangulation $\Delta$ of a half-dilation pillowcase, $\Phi_i(\Delta)$ is the marked triangle associated to the tetrahedral triangulation after flipping the edge that connects $p_i$ to $p_4$ and the edge that connects $p_{i+1 \mod 3}$ and $p_{i-1 \mod 3}$.
\begin{center}
\begin{tikzpicture}[scale = 0.45]

\filldraw[blue!20!white] (({2},0) -- (-1,{sqrt(3)}) -- (0,{2*sqrt(3)}) -- cycle;
\filldraw[yellow!20!white] (({2},0) -- (-1,{sqrt(3)}) -- (-2,0) -- cycle;
\draw[dotted] ({2},0) -- (0,{2*sqrt(3)});

\draw[latex-] (0,{2*sqrt(3)}) --({-2},0);
\draw[latex-] ({-2}, 0)-- ({2},0);
\draw[latex-] ({2}, 0)--  (0,{2*sqrt(3)}) ;

\draw[dashed] ({1},{sqrt(3)}) -- (-1,{sqrt(3)}) -- (0,0) -- (1,{sqrt(3)});
\draw (0,-0.5) node {$p_3$};
\draw ({1.5},{sqrt(3)+0.25}) node {$p_2$};

\draw[red] (2,0) -- (-1,{sqrt(3)});
\draw[red] (0,0) -- ({-1.5},{0.5*sqrt(3)});
\draw[red] ({1},{sqrt(3)}) -- (-0.5,{1.5*sqrt(3)});

\draw[->] (-1.25 - 0.1,{0.75*sqrt(3) + 0.1*(1/3)*sqrt(3)}) arc (240:60:0.5);

\draw (3.5, {sqrt(3)+0.5}) node {$\Phi_1$};
\draw[Bar->] (2.5, {sqrt(3)}) -- (4.5, {sqrt(3)});

\filldraw[blue!20!white] (({2+6},{2*sqrt(3)}) -- (4+6,{0}) -- (1+6,{sqrt(3)})  -- cycle;
\filldraw[yellow!20!white] (({-2+6},{2*sqrt(3)}) -- (1+6,{sqrt(3)}) -- (2+6,{2*sqrt(3)}) -- cycle;

\draw[-latex] (-2+6,{2*sqrt(3)}) --({4+6},0);
\draw[-latex] ({4+6}, 0)-- (2+6 ,{2*sqrt(3)}) ;
\draw[latex-] ({-2+6}, {2*sqrt(3)})--  ({2+6}, {2*sqrt(3)}) ;

\draw[dashed, red] ({1+6},{sqrt(3)}) -- (3+6,{sqrt(3)}) -- (0+6,{2*sqrt(3)}) -- cycle;
\draw ({3.5 + 6},{sqrt(3)+0.25}) node {$p_2$};
\draw ({7-0.5}, {sqrt(3)-0.25}) node {$p_1$};
\draw (0+6, {2*sqrt(3) + 0.5}) node {$p_3$};

\end{tikzpicture}\hspace{5pt} \begin{tikzpicture}[scale = 0.45]

\filldraw[blue!20!white] (({-2},0) -- (1,{sqrt(3)}) -- (0,{2*sqrt(3)}) -- cycle;
\filldraw[yellow!20!white] (({-2},0) -- (1,{sqrt(3)}) -- (2,0) -- cycle;
\draw[dotted] ({2},0) -- (0,{2*sqrt(3)});

\draw[latex-] (0,{2*sqrt(3)}) --({-2},0);
\draw[latex-] ({-2}, 0)-- ({2},0);
\draw[latex-] ({2}, 0)--  (0,{2*sqrt(3)}) ;

\draw[dashed] ({1},{sqrt(3)}) -- (-1,{sqrt(3)}) -- (0,0) -- (1,{sqrt(3)});
\draw (0,-0.5) node {$p_3$};
\draw ({-1.5},{sqrt(3)+0.25}) node {$p_1$};

\draw[red] (-2,0) -- (1,{sqrt(3)});
\draw[red] (0,0) -- ({1.5},{0.5*sqrt(3)});
\draw[red] ({-1},{sqrt(3)}) -- (0.5,{1.5*sqrt(3)});

\draw[->] (0.75 + 0.1,{1.25*sqrt(3) + 0.1*(1/3)*sqrt(3)}) arc (120:-60:0.5);

\draw (3.5, {sqrt(3)+0.5}) node {$\Phi_2$};
\draw[Bar->] (2.5, {sqrt(3)}) -- (4.5, {sqrt(3)});

\filldraw[blue!20!white] (({2+6},0) -- (4+6,{2*sqrt(3)}) -- (1+6,{sqrt(3)})  -- cycle;
\filldraw[yellow!20!white] (({-2+6},0) -- (1+6,{sqrt(3)}) -- (2+6,0) -- cycle;
\draw[dotted] ({2},0) -- (0,{2*sqrt(3)});

\draw[-latex] (4+6,{2*sqrt(3)}) --({-2+6},0);
\draw[-latex] ({2+6}, 0)-- (4+6,{2*sqrt(3)}) ;
\draw[-latex] ({-2+6}, 0)--  ({2+6}, 0) ;

\draw[dashed] ({1+6},{sqrt(3)}) -- (3+6,{sqrt(3)}) -- (0+6,0) -- cycle;
\draw (6,-0.5) node {$p_3$};
\draw ({3.5+6},{sqrt(3)-0.25}) node {$p_1$};
\draw ({0.5+6},{sqrt(3)+0.25}) node {$p_2$};

\end{tikzpicture}

\begin{center}
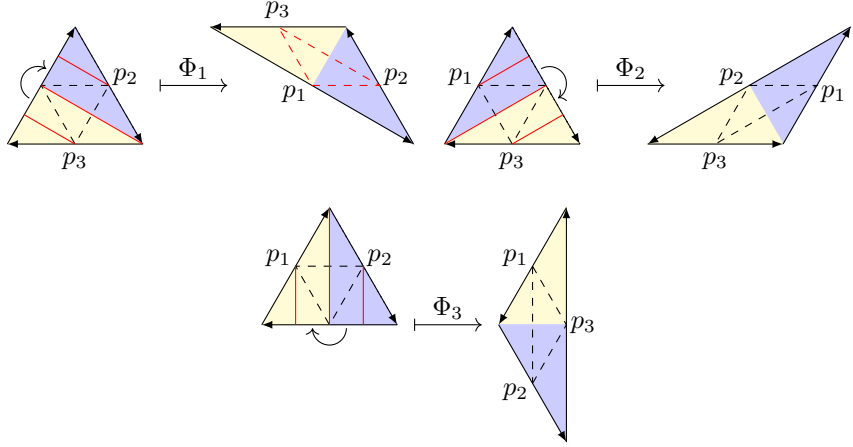

\begin{tikzpicture}[scale = 0.45]
\filldraw[blue!20!white] (({2},0) -- (0,{2*sqrt(3)}) -- (0,0) -- cycle;
\filldraw[yellow!20!white] (({-2},0) -- (0,{2*sqrt(3)}) -- (0,0) -- cycle;
\draw[dotted] ({2},0) -- (0,{2*sqrt(3)});

\draw[latex-] (0,{2*sqrt(3)}) --({-2},0);
\draw[latex-] ({-2}, 0)-- ({2},0);
\draw[latex-] ({2}, 0)--  (0,{2*sqrt(3)}) ;

\draw[dashed] ({1},{sqrt(3)}) -- (-1,{sqrt(3)}) -- (0,0) -- (1,{sqrt(3)});
\draw ({-1.5},{sqrt(3)+0.25}) node {$p_1$};
\draw ({1.5},{sqrt(3)+0.25}) node {$p_2$};

\draw[red] (0,0) -- (0,{2*sqrt(3)});
\draw[red] (1,0) -- ({1},{sqrt(3)});
\draw[red] (-1,0) -- ({-1},{sqrt(3)});

\draw[->] (0.5,-0.1) arc (0:-180:0.5);

\draw (3.5, 0.5) node {$\Phi_3$};
\draw[Bar->] (2.5, 0) -- (4.5, 0);

\filldraw[yellow!20!white] (({-2+7},0) -- (0+7,{2*sqrt(3)}) -- (0+7,0) -- cycle;
\filldraw[blue!20!white] (({-2+7},0) -- (0+7,{-2*sqrt(3)}) -- (0+7,0) -- cycle;
\draw[dashed] (0+6,{sqrt(3)}) -- (0+6,{-1*sqrt(3)});
\draw[dashed] (0+6,{sqrt(3)}) -- (7, 0);
\draw[dashed] (0+6,{-sqrt(3)}) -- (7, 0);

\draw[latex-] ({-2+7},0) -- (0+7,{2*sqrt(3)});
\draw[latex-] (0+7,{-2*sqrt(3)}) -- ({-2+7},0);
\draw[latex-] (0+7,{2*sqrt(3)}) -- (0+7,{-2*sqrt(3)});

\draw ({-1.5 + 7},{sqrt(3)+0.25}) node {$p_1$};
\draw ({-1.5 + 7},{-sqrt(3)-0.25}) node {$p_2$};
\draw ({0.5+7},{0}) node {$p_3$};

\end{tikzpicture}
\captionof{figure}{The three different flips}
\end{center}
\end{center}
We let the map $f_{i,\Delta}: \Delta \to \Phi_i(\Delta)$ be the actual cut-scale-paste operation.\nextpar

These three maps are important since, given a tetrahedral triangulation of a half-dilation pillowcase, the finite sequence of flips obtained from the flipping algorithm is actually a composition of these maps.

\begin{prop}
Let $T$ be a tetrahedral triangulation of a half-dilation pillowcase $X$, and $T_D$ the Delaunay triangulation of the same half-dilation pillowcase. Let $\set{T_k}_{k = 0}^N$ be the finite sequence of flips given by flipping algorithm gives a sequence of triangulated half-dilation surfaces where $T_0 = T$ and $T_N = T_D$. Then $N$ is even, and for each even $j$, there is some $i$ such that $T_{j+2} = \Phi_i(T_j)$.
\end{prop}

\begin{proof}
After the first Delaunay flip, a tetrahedral triangulation becomes a $\set{2,2,4,4}$-triangulation.
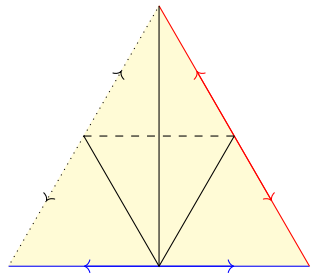
\begin{center}
\begin{tikzpicture}

\filldraw[yellow!20!white] (({2},0) -- (0,{2*sqrt(3)}) -- ({-2},0) -- ({2},0) -- cycle;

\draw[blue] ({2},0) -- (-2,0);
\draw[blue, ->] (0,0) -- (1,0);
\draw[blue, ->] (0,0) -- (-1,0);

\draw[red, ->] ({1},{sqrt(3)}) -- (1.5, {0.5*sqrt(3)});
\draw[red, ->] ({1},{sqrt(3)}) -- (0.5, {1.5*sqrt(3)});
\draw[red] ({2},0) -- (0,{2*sqrt(3)});

\draw[dotted, ->] ({-1},{sqrt(3)}) -- (-1.5, {0.5*sqrt(3)});
\draw[dotted, ->] ({-1},{sqrt(3)}) -- (-0.5, {1.5*sqrt(3)});
\draw[dotted] (-1.5, {0.5*sqrt(3)}) -- (-2,0);
\draw[dotted] (-.5, {1.5*sqrt(3)}) -- (0,{2*sqrt(3)});

\draw[-] (-1,{sqrt(3)}) -- (0,0) -- (1,{sqrt(3)});
\draw[-] (0,0) -- (0,{2*sqrt(3)});

\draw[dashed] ({1},{sqrt(3)}) -- (-1,{sqrt(3)});

\end{tikzpicture}
\end{center}

Note that the only edge that is potentially not locally Delaunay is the one blue edge. Indeed, the angles opposite the red edge (or dotted edge) form a corner of a triangle, so its sum must be less than $\pi$, and the angles opposite the two internal edges other than the edge that was just created form two angles of a triangle, which again, must be less than $\pi$. We conclude that the $(2n-1)$-th and $2n$-th Delaunay flips are on disjoint pairs of triangles.\nextpar

Furthermore, there must be an even number of flips, since every odd number of flips takes a tetrahedral triangulation to a $\set{2,2,4,4}$-triangulation. We conclude that there are an even number of Delaunay flips, and if we pair the $2n$-th and $(2n+1)$-th flip, we get some $\Phi_i$.\nextpar

Now, we further show that $\frac{N}{2}$ is even. Let $\set{T_{2k}}_{k = 0}^{\frac{N}{2}}$ be the subsequence of tetrahedrally triangulated half-dilation surfaces. Then $T_{2(k + 1)} = \Phi_{i_k}(T_{2k})$. Because each $\Phi_{i_k}$ reverses the orientation of the dilation ratios, there must be an even number of $\Phi_{i_k}$'s. 
\end{proof}

We now consider the linear maps associated to each $\Phi_i$. In particular, given a tetrahedrally triangulated pillowcase $X$, let $\Delta \in \Delta^*$ be a marked triangle associated to $X$. Then for each $i = 1, 2, 3$ there is a unique affine map $\phi_{i,\Delta} = Ax + b : \R^2 \to \R^2$ where $A \in \GL(2, \R)$ and $b \in \R^2$ such that
\[
\phi_{i,\Delta}(\Delta) = \Phi_i(\Delta).
\]
Letting $\cal{T}^{-1}(\Delta) = \big( (v_1, \lambda_1), (v_2, \lambda_2), (v_3, \lambda_3) \big)$ and $\cal{T}^{-1}(\Phi_i) = \big( (v_1', \lambda_1), (v_2', \lambda_2), (v_3', \lambda_3) \big)$, the linear map $\phi_{i,\Delta}$ must be the one that takes $v_1 \mapsto v_1'$ and $v_2 \mapsto v_2'$. \nextpar

We can actually define these $\phi_{i, \Delta}$ explicitly by determining the what vectors are after making the flip $\Phi_i$. Let $\cal{T}^{-1}(\Delta) = \big( (v_1, \lambda_1), (v_2, \lambda_2), (v_3, \lambda_3) \big)$
\begin{align*}
\cal{T}^{-1}(\Phi_1(\Delta)) &= \big((-(\lambda_1+1)v_3-v_1, \lambda_1), (-v_2, \lambda_2), (\lambda_1v_3, \lambda_3)\big)\\
\cal{T}^{-1}(\Phi_2(\Delta)) &= \big((\lambda_2v_1, \lambda_1), (-(\lambda_2+1)v_1-v_2, \lambda_2), (-v_3, \lambda_3)\big)\\
\cal{T}^{-1}(\Phi_3(\Delta)) &= \big((-v_1, \lambda_1), (\lambda_3v_2, \lambda_2), (-(\lambda_3+1)v_2 - v_3, \lambda_3)\big)
\end{align*}
Using the vectors $v_1$ and $v_2$ of the tetrahedral triangulation $\Delta$, we can write a matrix $\cal{M}(\phi_{i, \Delta})$ down for each map $\phi_{i, \Delta}$. Let $P = \begin{bmatrix} v_1 & v_2\end{bmatrix}$.
\begin{align*}
\cal{M}(\phi_{1, \Delta}) &= P \circ \begin{bmatrix}
\lambda_1 & 0\\
\lambda_1+1  & -1
\end{bmatrix} \circ P^{-1} , \\
\cal{M}(\phi_{2, \Delta}) &=  P \circ \begin{bmatrix}
\lambda_2 & 0\\
 -(\lambda_2+1)  & -1
\end{bmatrix} \circ P^{-1}, \\
\cal{M}(\phi_{3, \Delta}) &= P \circ \begin{bmatrix}
-1 & 0\\
0  & \lambda_3
\end{bmatrix} \circ P^{-1}
\end{align*}
Using these maps, we can construct some affine automorphisms. Given a tetrahedrally triangulated pillowcase $X$, let $\Delta \in \Delta^*$ be a marked triangle associated to $X$. Restricting the map $\phi_{i,\Delta}$ to $\Delta$, consider the composition of the maps,
\begin{center}
\begin{tikzcd}
\Delta \ar[r, "\phi_{i, \Delta}"] & \Phi_i(\Delta) \ar[r, "\phi_{j, \Phi_i(\Delta)}"] & \Phi_j\Phi_i(\Delta) \ar[r, "f_{j, \Phi_i(\Delta)}^{-1}"] & \Phi_i(\Delta) \ar[r, "f_{i, \Delta}^{-1}"] & \Delta
\end{tikzcd}
\end{center}

To make sure this is an affine automorphism, we take its derivative. The map $f_{i, \Delta}^{-1}\circ f_{j, \Phi_i(\Delta)}^{-1}$ corresponds to a cut-paste operation, which is a local half-dilation. We now write down the linear maps $\phi_{j, \Phi_i(\Delta)} \circ \phi_{i, \Delta}$ by looking at $\cal{T}^{-1}(\Phi_j\Phi_i(\Delta))$.
\begin{align*}
\cal{T}(\Phi_3\Phi_2(\Delta)) &= \big(-\lambda_2v_1, \lambda_1), (-\lambda_3(\lambda_2+1)v_3 - \lambda_3v_1 , \lambda_2), ((\lambda_2 + 1)(\lambda_3 + 1)v_1 + (\lambda+1)v_2 + v_3, \lambda_3)\big)\\
\cal{T}(\Phi_2\Phi_1(\Delta)) &= \big((-\lambda_2(\lambda_1+1)v_3 - \lambda_2v_1, \lambda_1), ((\lambda_1 + 1)(\lambda_2 + 1)v_3 + (\lambda+1)v_1 + v_2 , \lambda_2), (-\lambda_1v_3, \lambda_3)\big)\\
\cal{T}(\Phi_1\Phi_3(\Delta)) &= \big(((\lambda_3 + 1)(\lambda_1 + 1)v_2 + (\lambda+1)v_3 + v_1 , \lambda_1), (-\lambda_3v_2, \lambda_2), (-\lambda_1(\lambda_3+1)v_2 - \lambda_1v_3, \lambda_3)\big)
\end{align*}
From this, denoting $\phi_{j, \Phi_i(\Delta)}\phi_{i, \Delta}$ as $\psi_{j, i, \Delta}$ and letting $P_{i,j} = \begin{bmatrix} v_i & v_j \end{bmatrix}$ we can conclude that,
\begin{align*}
\cal{M}(\psi_{3, 2, \Delta}) &= P\circ \begin{bmatrix}
-\lambda_2 & -\lambda_3(\lambda_2+1) \\
0 & -\lambda_3
\end{bmatrix} \circ P^{-1},\\
\cal{M}(\psi_{2, 1, \Delta}) &= P\circ \begin{bmatrix}
 \lambda_1\lambda_2 & \lambda_2(\lambda_1 + 1)\\
-\lambda_1(\lambda_2+1) & -\lambda_1\lambda_2-\lambda_1-\lambda_2
\end{bmatrix} \circ P^{-1},\\
\cal{M}(\psi_{1, 3, \Delta}) &= P\circ \begin{bmatrix}
\lambda_3\lambda_1 & 0  \\
\lambda_3(\lambda_1 + 1) & 1
\end{bmatrix} \circ P^{-1}
\end{align*}
When we derive $f^{-1}_{i,\Delta}\circ f_{j, \Phi_i(\Delta)}^{-1} \circ  \phi_{j, \Phi_i(\Delta)} \circ \phi_{i,\Delta}$, we must use the chain rule. Since the derivative a local half-dilation is a half-dilation, and the derivative of a linear map is itself, we get that the derivative is
\[
\begin{bmatrix}
0 &-1\\
-1& 0
\end{bmatrix}^n \cdot \begin{bmatrix}
\lambda & 0\\
0 & \lambda
\end{bmatrix} \cdot \cal{M}^{v_1}_{v_2}(\psi_{j, i, \Delta}) .
\]
Taking its image under the quotient map, we get that $\bof{\cal{M}^{v_1}_{v_2}(\psi_{j, i, \Delta})}$ 

\begin{prop}
$\PSL(X)$ is generated by $\set{\bof{\cal{M}^{v_1}_{v_2}(\psi_{j, i, \Delta_D})}}$.
\end{prop}
\begin{proof}
Given $X \in \mathcal{HD}(S_{0,4})$, let $\Delta_D$ be its Delaunay triangulation and $[A] \in \PSL(X)$. It suffices to show that there is a sequence of flips $\Phi_i\Phi_j$ taking $\cal{T}^{-1}(A\cdot \cal{T}(\Delta_D))$ to $\Delta_D$.\nextpar

If $[A]$ is in the Veech group, then $\cal{T}^{-1}(A\cdot \cal{T}(\Delta_D))$ must also be a tetrahedral triangulation of $X$. By Theorem 2.11 there is a finite sequence of Delaunay flips of locally non-Delaunay edges taking any geodesic triangulation to a Delaunay triangulation. Call this sequence of Delaunay flips $\set{f_n}$. By propositions 5.4 and 5.5, we can conclude that there is a finite sequence of $\Phi_i\Phi_j$'s taking $A\cdot X$ to $X$. We conclude that $A$ is some composition of $\cal{M}^{v_1}_{v_2}(\phi_{j, i, \Delta_D})$'s.
\end{proof}

\section{Connection to Hyperbolic Structures on $S_{0,3}$}
Recall that we can classify matrices in $\SL(2,\R)$ by elliptic, parabolic, and hyperbolic according to its trace. We calculate this for the matrices $\cal{M}^{v_1}_{v_2}(\phi_{j, i, \Delta_D}) \in \GL(2,\R)$ by normalizing by multiplying the matrix $A$ with $\sqrt{\tfrac{1}{\det(A)}}$ and taking the trace to get $\tfrac{\tr(A)}{\sqrt{\det(A)}}$. For convenience, we take the square of this conjugacy invariant to calculate the following:

\begin{align*}
\frac{\tr^2}{\det}\of{\cal{M}^{v_1}_{v_2}(\phi_{3, 2, \Delta_D})} &= \frac{\tr^2}{\det}\begin{bmatrix}
-\lambda_2 & -\lambda_3(\lambda_2+1) \\
0 & -\lambda_3
\end{bmatrix} &= \lambda_3\lambda_2 + 2 + \frac{1}{\lambda_3\lambda_2}\\
\frac{\tr^2}{\det}\of{\cal{M}^{v_1}_{v_2}(\phi_{2, 1, \Delta_D})} &= \frac{\tr^2}{\det}\begin{bmatrix}
 \lambda_1\lambda_2 & \lambda_2(\lambda_1 + 1)\\
-\lambda_1(\lambda_2+1) & -\lambda_1\lambda_2-\lambda_1-\lambda_2
\end{bmatrix} &= \lambda_2\lambda_1 + 2 + \frac{1}{\lambda_2\lambda_1}\\
\frac{\tr^2}{\det}\of{\cal{M}^{v_1}_{v_2}(\phi_{1, 3, \Delta_D})}&= \frac{\tr^2}{\det}\begin{bmatrix}
\lambda_3\lambda_1 & 0  \\
\lambda_3(\lambda_1 + 1) & 1
\end{bmatrix} &= \lambda_1\lambda_3 + 2 + \frac{1}{\lambda_1\lambda_3}
\end{align*}

The trace gives us information about the hyperbolic structure on $S_{0,3}$ derived from $\mathbb{H}^2/\PSL(X)$. In particular, it gives us the lengths of curves around each of the singularities through the following relationship:
\[
\tr(A) = 2\cosh(l)
\]
where $l$ is the length of the curve associated to $A$. This allows us to define a map to Fenchel-Nielsen coordinates to the space of complete hyperbolic structures on $S_{0,3}$. If there is a triple of dilation ratios whose normalized square trace is $(l_a, l_b, l_c) \in [4, \infty)^3$ we can conclude that every hyperbolic structure on $S_{0,3}$ is represented by a quotient of $\mathbb{H}^2$ by the Veech group of a half-dilation structure on $S_{0,4}$.

\begin{prop}
The map $\R_{>0}^3 \to [4, \infty)^3$ defined by 
\[
(\lambda_1, \lambda_2, \lambda_3) \mapsto \left(\lambda_3\lambda_2 + 2 + \frac{1}{\lambda_3\lambda_2}, \lambda_2\lambda_1 + 2 + \frac{1}{\lambda_2\lambda_1}, \lambda_1\lambda_3 + 2 + \frac{1}{\lambda_1\lambda_3}\right)
\]
is surjective.
\end{prop}
\begin{proof}
Let $a, b, c \in [4, \infty)$. We solve the following system of equations for $\lambda_1, \lambda_2, \lambda_3$:
\[
\begin{cases}
a &= \lambda_3\lambda_2 + 2 + \frac{1}{\lambda_3\lambda_2}\\
b &= \lambda_2\lambda_1 + 2 + \frac{1}{\lambda_2\lambda_1}\\
c &= \lambda_1\lambda_3 + 2 + \frac{1}{\lambda_1\lambda_3}
\end{cases}
\Leftrightarrow
\begin{cases}
a-2 &= \lambda_3\lambda_2 + \frac{1}{\lambda_3\lambda_2}\\
b-2 &= \lambda_2\lambda_1 + \frac{1}{\lambda_2\lambda_1}\\
c-2&= \lambda_1\lambda_3 + \frac{1}{\lambda_1\lambda_3}
\end{cases}
\]
since $y = x + \tfrac{1}{x}$ has the solutions, $\frac{y \pm \sqrt{y^2 - 4}}{2}$. we can simplify the system of equations to,
\[
\begin{cases}
\lambda_3\lambda_2 &= \frac{(a-2) \pm \sqrt{a^2 - 4a}}{2}\\
\lambda_2\lambda_1 &= \frac{(b-2) \pm \sqrt{b^2 - 4b}}{2}\\
\lambda_1\lambda_3 &= \frac{(c-2) \pm \sqrt{c^2 - 4c}}{2}
\end{cases}
\]
letting $s_{\pm}(x) = \frac{(x-2) \pm \sqrt{x^2 - 4x}}{2}$, we conclude that
\[
\begin{cases}
\lambda_1 &= \sqrt{\frac{s_{\pm}(a)s_{\pm}(c)}{s_{\pm}(b)}}\\
\lambda_2 &= \sqrt{\frac{s_{\pm}(a)s_{\pm}(b)}{s_{\pm}(c)}}\\
\lambda_3 &= \sqrt{\frac{s_{\pm}(b)s_{\pm}(c)}{s_{\pm}(a)}}
\end{cases}
\]
so the map is surjective.
\end{proof}
By choice of $\pm$ for each $a, b, c$, there is at most 8 elements that map to each triple in $[4, \infty)^3$. We explore what the relationship is between these 8 triples of dilation ratio. First, note that $s_+(x)s_-(x) = 1$, so $s_+(x)$ and $s_-(x)$. are reciprocals of each other. Let 
\[
\begin{cases}
\lambda_1 &= \sqrt{\frac{s_{+}(a)s_{+}(c)}{s_{+}(b)}}\\
\lambda_2 &= \sqrt{\frac{s_{+}(a)s_{+}(b)}{s_{+}(c)}}\\
\lambda_3 &= \sqrt{\frac{s_{+}(b)s_{+}(c)}{s_{+}(a)}}
\end{cases}
\]
Then,
\[
\frac{1}{\lambda_4} = \lambda_1\lambda_2\lambda_3 =  \sqrt{\frac{s_{+}(a)s_{+}(c)}{s_{+}(b)}} \sqrt{\frac{s_{+}(a)s_{+}(b)}{s_{+}(c)}} \sqrt{\frac{s_{+}(b)s_{+}(c)}{s_{+}(a)}} = \sqrt{s_+(a)s_+(b)s_+(c)}.
\]
If we consider the solution where $s_+(a)$ is replaced with $s_-(a)$, we get
\[
\begin{cases}
\lambda_1' &= \sqrt{\frac{s_{-}(a)s_{+}(c)}{s_{+}(b)}} = \sqrt{\frac{s_+(c)}{s_+(a)s_+(b)}} = \frac{1}{\lambda_2} \\
\lambda_2' &= \sqrt{\frac{s_{-}(a)s_{+}(b)}{s_{+}(c)}} = \sqrt{\frac{s_+(b)}{s_+(a)s_+(c)}} = \frac{1}{\lambda_1}\\
\lambda_3' &= \sqrt{\frac{s_{+}(b)s_{+}(c)}{s_{-}(a)}} = \sqrt{s_+(a)s_+(b)s_+(c)} = \frac{1}{\lambda_4}.
\end{cases}
\]
After doing very similar calculations, we get that the solutions are marked triangulations with the following dilation ratios:
\begin{align*}
(\lambda _1, \lambda_2, \lambda_3) && (\tfrac{1}{\lambda_1}, \tfrac{1}{\lambda_2}, \tfrac{1}{\lambda_3})\\
(\lambda _2, \lambda_1, \lambda_4) && (\tfrac{1}{\lambda_2}, \tfrac{1}{\lambda_1}, \tfrac{1}{\lambda_4})\\
(\lambda _4, \lambda_3, \lambda_2) && (\tfrac{1}{\lambda_4}, \tfrac{1}{\lambda_3}, \tfrac{1}{\lambda_2})\\
(\lambda _3, \lambda_4, \lambda_1) && (\tfrac{1}{\lambda_3}, \tfrac{1}{\lambda_4}, \tfrac{1}{\lambda_1})
\end{align*}
However, marked triangulation with reciprocal dilation ratios differ by the linear map $\begin{bmatrix} 0 & 1\\ 1& 0\end{bmatrix}$. Therefore, these triples of dilation ratios of marked triangulations correspond to an unmarked triangulation with 3 dilation ratios. 
\begin{cor}
Given $P$, a triangulable half-dilation structure on the four-punctured sphere, it's Veech Group $\PSL(P) < \PSL(2, \R)$ is such that $\mathbb{H}^2/\PSL(P)$ is a complete hyperbolic structure on the three-punctured sphere. Furthermore, every hyperbolic structure on the three-punctured sphere is represented.
\end{cor}

\subsection*{Acknowledgements}
This work was partially supported by the Dr. Barnett and Jean Hollander Rich Mathematics Scholarship.
The author is grateful for W. Patrick Hooper for his guidance and thoughts.

\bibliographystyle{alpha}
\bibliography{ref}

\end{document}